\documentclass[11pt]{amsart}

\usepackage{amsmath, amssymb, amscd}
\usepackage{eucal}
\usepackage{mathrsfs}
\usepackage[frame,cmtip,arrow,matrix,line,graph,curve]{xy}
\usepackage{graphicx}

\newcommand{\rr}{\mathbb R}
\newcommand{\zz}{\mathbb Z}

\newtheorem{prop}{Proposition}[section]
\newtheorem{thm}[prop]{Theorem}
\newtheorem{lem}[prop]{Lemma}
\newtheorem{cor}[prop]{Corollary}

\theoremstyle{remark}
  \newtheorem{rk}[prop]{Remark}
  
\theoremstyle{definition}
 \newtheorem{exam}[prop]{Example}
 \newtheorem{defn}[prop]{Definition}

\numberwithin{equation}{section}

\begin{document}
\title[On Homoclinic points, Recurrences and Chain recurrences of
volume-preserving diffeomorphisms without genericity] {On Homoclinic
points, Recurrences and Chain recurrences of volume-preserving
diffeomorphisms without genericity}
\author{Jaeyoo Choy, Hahng-Yun Chu$^{\ast}$ and Min Kyu Kim}
\address{Dept. of Mathematics, Kyungpook National University, Sankyuk-dong, Buk-gu,
Daegu 702-701, Republic of Korea}\email{choy@knu.ac.kr}
\address{School of Mathematics, Korea Institute for Advanced Study, Cheongnyangni 2-dong, Dongdaemun-gu,
Seoul 130-722, Republic of Korea}\email{hychu@kias.re.kr}
\address{Department of Mathematics Education,
Gyeongin National University of Education, Gyesan-dong,
Gyeyang-gu, Incheon, 407-753, Republic of Korea} \email{mkkim@kias.re.kr}

\thanks{\it $\ast$ Corresponding author.}
\thanks{The first and second authors were supported by Korea Research Foundation
Grant(KRF-2008-331-C00015)
}

\subjclass[2000]{primary 37C50; secondary 37C29, 37D05}
\keywords{recurrent point, homoclinic
point, hyperbolic, shadowable, chain recurrence, Lagrange-stable}

\begin{abstract}

Let $M$ be a manifold with a volume form $\omega$ and $f : M
\rightarrow M$ be a diffeomorphism of class $\mathcal{C}^1$ that
preserves $\omega$. In this paper, we do \textit{not} assume $f$ is
$\mathcal{C}^1$-generic. We have two main themes in the paper: (1)
the chain recurrence; (2) relations among recurrence points,
homoclinic points, shadowability and hyperbolicity. For (1) (without
assuming $M$ is compact), we have the theorem: if $f$ is Lagrange
stable, then $M$ is a chain recurrent set. If $M$ is compact, then
the Lagrange-stability is automatic. For (2) (assuming the
compactness of $M$), we prove some various implications among
notions, such as: (i) the $\mathcal{C}^1$-stable shadowability
equals to the hyperbolicity of $M$; (ii) if a point $p\in M$ has a
recurrence point in the unstable manifold $W^u (p, f)$ and there is
no homoclinic point of $p,$ then $f$ is nonshadowable; (iii) if $f$
has the shadowing property and $p$ has a recurrence point in $W^u
(p, f),$ then the recurrent point is in the limit set of homoclinic
points of $p$.
\end{abstract}

\maketitle

\thispagestyle{empty} \markboth{Jaeyoo Choy, Hahng-Yun Chu and Min
Kyu Kim} {On Homoclinic points, Recurrences and Chain recurrences
without genericity}


\section{Introduction}\label{sec: intro}

An important part of the study of Dynamical system is to understand
the structures of the stable manifolds (and the unstable manifolds),
attractors and chain recurrent sets of the dynamical system. In
order to grasp the geometric structure in Dynamical system, the
above concepts have played an important role in the field bearing
the useful properties.

Our purpose of this paper is to elucidate \textit{the relation
between the set of homoclinic points and the set of  recurrence
points}, and to know \textit{the existence of the chain recurrence
points}. In the both themes, we want to maintain the
volume-preserving (or symplectic) dynamics, but to drop the
genericity assumption. Since the symplectic diffeomorphisms are
automatically volume-preserving, our results in the paper is
applicable to the symplectic dynamics as well.

In this paper, we follow the two basic definitions of attractors and
chain recurrences due to Conley \cite{Co}. Note that even before
Conley's definition, other definitions of attractor can be found in
several papers (see, for reference, \cite{Mi}) while the equivalence
of those definitions partly remain as conjecture.

By a theorem of Moriyasu \cite{Mo}, one can see that if the chain
recurrent set $CR(f)$ is $\mathcal{C}^1$-stably shadowable, then it
is hyperbolic (see \cite{WGW}). In \cite{NN}, Nakayama and Noda
proved that for a minimal flow on a closed orientable $3$-manifold,
a chain recurrent set of the induced projective flow is the whole
projectivized bundle. Recently, it is proved in \cite{WGW} that if a
chain component is $\mathcal{C}^1$-stably shadowable, then the chain
component is hyperbolic. Other results for  the
$\mathcal{C}^1$-stable shadowing properties can be found in
\cite{LMS} and \cite{Sa}.

We obtain our first result of this article as follows.\\

\noindent \textbf{Theorem \ref{thm:global hyp}.} \textit{Let $M$
be a compact manifold with a volume form $\omega$ and $f$ be a
volume-preserving diffeomorphism on $M$. Assume that every chain
component has a hyperbolic periodic point. If every chain
component is $\mathcal{C}^{1}$-stably shadowable, then $M$ is
hyperbolic. Furthermore, the converse is also true.}\\

By the diffeomorphisms throughout this paper, we always mean the $\mathcal C^1$-diffeomorphisms.

From Conley's theorem, we get easily the Chain recurrence theorem
on the compact manifolds with volume-preserving diffeomorphisms,
i.e., $M$ is chain recurrent for $f$. Furthermore, in fact, $M$ is
nonwandering for $f$ by the Nonwandering theorem. Howerer,
unfortunately, the same statement for the chain recurrence does
not hold for the non-compact manifolds with volume-preserving
diffeomorphisms. In the non-compact case, we may consider a
canonical notion, Lagrange stability, and using this concept, we
prove  the following
theorem.\\

\noindent \textbf{Theorem \ref{th: LS CR}.} \textit{Let $M$ be a
manifold with a volume form $\omega$, and $f$ be a Lagrange-stable
volume-preserving diffeomorphism on $M$. Then, $M$ is strongly chain
recurrent for $f$, i.e., every point of $M$ is strongly chain
recurrent with respect to $f$.}\\

\noindent The above theorem follows from a stronger claim,
Proposition \ref{prop: measure 0}. The proposition says that if the
assumptions in Theorem \ref{th: LS CR} without the
Lagrange-stability assumption hold, then almost everywhere in $U-A$
should have the unbounded orbit, where $A$ is an attractor and $U$ is an
attractor block of $A$. I.e., the set of points of $U-A$ with
bounded orbits is of measure 0.

Once focusing on dynamics on symplectic geometry, so called
Hamiltonian dynamics, one realizes that the study of homoclinic
points sits in a weighty place, which dates back to a century
ago(see \cite{BGW}). More generally, Poincar$\acute{e}$ raised a
question:   \textit{do transverse homoclinic points occur
generically?} It is interpreted as whether transverse homoclinic
points are dense in both stable and unstable manifolds. He
conjectured that the answer is yes. Xia \cite{Xi1} answered
positively for the Poincar$\acute{e}$ question on arbitrary compact
manifolds with a generic volume-preserving diffeormorphism.
His
results in \cite{Xi1} generalize Takens' result \cite{Ta}: in the generic
volume-preserving dynamics on a compact manifold (it is called a
\textit{discrete conservative system} in \cite{Ta}; Takens also
treats Hamiltonian flows, a \textit{continuous conservative
system}), for a hyperbolic periodic point,
the set of all homoclinic points of the hyperbolic periodic point are dense in both the stable
and unstable manifolds.

In \cite{Xi2}, for the compact surface case (with volume-preserving
diffeormorphisms), Xia shows even that the union of the stable
manifolds(unstable manifolds) of hyperbolic periodic points are dense in the whole manifold.

Furthermore, it is proved that the existence of a residual subset $\mathcal{R}$ in the set of volume-preserving
diffeomorphisms such that for any $f \in \mathcal{R}$ and any hyperbolic periodic point of $f$,
the set of corresponding homoclinic points is dense in $M$ (\cite{SX}).

There is a parallel question with the above discussion on homoclinic
points: \textit{how generically do recurrence points occur?}
It turns out, due to Oliveira \cite{Ol2}, that the
answer is that for a generic $\mathcal{C}^1$-diffeomorphism
preserving $\omega$, the positive(negative) recurrence points form a dense subset in
the stable manifolds(unstable manifolds) of hyperbolic period points.

Lastly, in this paper, we prove the following theorem.\\

\noindent \textbf{Theorem \ref{th: rec nonsh}.} \textit{Let $M$ be a
compact manifold with a volume form $\omega$ and $f : M \rightarrow
M$ be a volume-preserving diffeomorphism. If $p$ has a recurrence point in $W^u (p, f)$ and there is
no homoclinic point of $p,$ then $f$ is nonshadowable.}\\

In \cite{BGW}, for $\mathcal{C}^1$-generic diffeomorphisms, it is
proven that every chain recurrent class that has a partially
hyperbolic splitting $E^s\oplus E^c\oplus E^u$ with $\dim E^c=1$
\textit{either} is an isolated hyperbolic periodic orbit,
\textit{or} is accumulated by non-trivial homoclinic classes. For
the recent investigations of the theory of volume-preserving or
symplectic dynamics, see, for instance, \cite{SX, XZ}.\\

\textit{Acknowledgement.} The second author is grateful to Dr. X.
Wen who explained him a corollary of Moriyasu's theorem in
\cite{WGW}.


\section{Definitions and Notations}\label{sec: def}

We fix the notations and definitions used throughout the paper.

Let $M$ be an $n$-dimensional differentiable manifold with a metric
$d$, and $f: M \rightarrow M$ be a $\mathcal{C}^1$-diffeomorphism. A
\textit{volume form} $\omega$ on $M$ is a nowhere vanishing $n$-form
on $M$. A symplectic form $\omega$ on $M$ is a nowhere degenerate
2-form on $M$. Here, the non-degeneracy of $\omega$ is the same as
its $(n/2)$-times wedge product $\omega^{\frac n2}=
\omega\wedge\cdots\wedge\omega$ defines a volume form on $M$. Thus,
when we say a symplectic form, $n$ is assumed to be even.
Integration along the subsets of $M$ defines a Lebesgue measure $m$.
Indeed, by the para-compactness of $M$, locally $m$ is written as a
product of a $\mathcal{C}^1$-function and the standard Lebesgue
measure on $\rr^n$ (via the $\mathcal{C}^1$-transition). This
clarifies a Lebesgue measurable subset of $M$,  a countable union of
Lebesgue measurable subsets of $\rr^n$ (via the
$\mathcal{C}^1$-transition). Thanks to the well-known theory of
Lebesgue measures and Borel measures, one guarantees any compact
subset of $M$ is Lebesgue measurable and is of finite measure. By
the compactness, the closed balls (with finite radii) are of finite
measure, as well.

If one says $f$ preserves $\omega$, this means $f^*\omega=\omega$.
When $\omega$ is a symplectic form, the $\omega$-preservation
implies the volume-preservation. The volume-preservation of $f$
amounts to  the measure-preservation. In the case, for a Lebesgue
measurable subset $N\subset M$, we have $m(N)=m(f(N))$.

\smallskip

For $p \in M$, the \textit{stable} and \textit{unstable manifolds}
are defined by
$$
W^s (p,f) = \{ u \in X : d(f^n(u), f^n(p)) \rightarrow 0 \ as \
n\rightarrow \infty \},
$$
$$
W^{u}(p, f) = \{ u \in X : d(f^{-n}(u), f^{-n}(p)) \rightarrow 0  \
as \ n\rightarrow \infty \},
$$ respectively.
\smallskip

\begin{defn}
A sequence $(x_i)_{i \in \zz}$ in $M$ is a
\textit{$\delta$-pseudo-orbit} of $f$ if for all $i \in \zz$,
$$d(x_{i+1},f(x_i))<\delta. $$

Given $\epsilon >0,$ a pseudo-orbit $(x_i)$ is
\textit{$\epsilon$-shadowed} by an actual orbit $(f^i (x))_{i \in
\zz}$ of $f$ if for all $i \in \zz,$ $$d(x_i,f^i (x))<\epsilon.$$

The diffeomorphism $f$ has the \textit{shadowing property} if for
every $\epsilon >0$, there exists some $\delta>0$ such that every
$\delta$-pseudo-orbit is $\epsilon$-shadowed by some actual orbit of
$f.$ In contrast, a diffeomorphism $f$ is \textit{nonshadowable} if
there is some $\epsilon
>0$ such that for every $\delta >0$, there is some $\delta$-pseudo-orbit which
is not $\epsilon$-shadowed by any actual orbit in $M.$ \end{defn}

\smallskip

\begin{defn}
A linear map of $\mathbb{R}^{n}$ is called \textit{hyperbolic} if it
admits $n$ eigenvectors (counted with multiplicities) and all the
eigenvalues have absolute values different from one. A fixed point
$p$ of a differentiable map $f$ is called \textit{hyperbolic} if the
tangent map $Df_{p}$ is hyperbolic. A periodic point $p$ of $f$ with
period $n$ is called \textit{hyperbolic periodic} for $f$ if
$D(f^{n})_{p} : T_{p}M \rightarrow T_{p}M$ is hyperbolic. We call
the orbit of such $p$ a \textit{hyperbolic periodic orbit.}

\end{defn}

\smallskip

\begin{defn} A point $q$ is called \textit{homoclinic} if it is in the
intersection of the stable manifold $W^s(p,f)$ and the unstable
manifold $W^u(p,f)$ of the same hyperbolic periodic point $p$ and
$p\neq q$. \end{defn}

\smallskip

\begin{defn}
A compact $f$-invariant set $\Lambda$ is called \textit{hyperbolic}
if the tangent bundle $T_{\Lambda}M$ has a continuous $Df$-invariant
splitting $E \oplus F$ and there exist constants $C > 0$, $0<
\lambda <1$ such that $\|Df^{n}|_{E(x)} \| \leq C\lambda^{n}$ and
$\|Df^{-n}|_{F(f^{n}(x))} \| \leq C\lambda^{n}$ for all $x \in
\Lambda$ and $n \geq 0$.
\end{defn}

\smallskip

The following definitions do not require the base space $X$ is a
manifold.

\begin{defn}\label{def: limit}
Let $f$ be a homeomorphism of a compact metric space $X$ and $p \in
X$. We define the $\textit{omega-limit set}$ of $p$, denoted as
$\omega(p,f)$,
 by
$$\{ q \in M : q=\lim_{j\rightarrow\infty} f^{n_j}(p)
\mbox{ \ for some sequence of integer}\ n_{j}\rightarrow \infty
\}.$$
\end{defn}

The set $\omega(p,f)$ is closed, nonempty and $f$-invariant, i.e.,
$f(\omega(p,f))= \omega(p,f)$.

\smallskip

\begin{defn} A point $q \in X$ is \textit{(positively) recurrent} or
\textit{omega-recurrent} with respect to $f$ if $q \in
\omega(q,f)$.\end{defn}


\section{$\mathcal{C}^{1}$-stably shadowing theorem}\label{sec: chain}

In this section, on a compact manifold $M$ with a volume-preserving
diffeomorphism, we study chain recurrence sets and attractors. And,
we prove that if $M$ is
$\mathcal{C}^{1}$-stably shadowable, then $M$ is hyperbolic. I.e.,
we can state that the notion of hyperbolicity is expressed in terms of
the shadowing property with the conditions for robustness in
volume-preserving diffeomorphisms. Therefore, the notions of
hyperbolicity and $\mathcal{C}^{1}$-stable shadowablity  coincide.

Firstly, we consider a Chain Recurrence Theorem and Nonwandering
Theorem on compact manifolds and then show that on the connected
compact manifolds $M$ with volume-preserving diffeomorphisms, $M$
is chain transitive (thanks to a theorem of Bonatti and et al).

Let $M$ be a compact manifold and $\omega$ be a volume form. Let
$Diff^{1}_{\omega}(M)$ be the set of volume-preserving
$\mathcal{C}^{1}$-diffeomorphisms of $M$. Bonatti and et
al.\cite{BC} proved that for $\mathcal{C}^1$-generic
diffeomorphisms, the limit set, the non-wandering set and the chain
recurrence set coincide, and that they are the closure of the set of
hyperbolic periodic points. Also they showed that
(assuming that the manifold is connected) there is a residual subset of $Diff^{1}_{\omega}(M)$ consisting
of transitive diffeomorphisms.

The definitions below need not assume the base space is a manifold.
Let $(X, d)$ be a metric space and $f:X \rightarrow X$ be a
continuous map.

\smallskip

\begin{defn}Let $\varepsilon >0$. Then a nonempty open subset $U$ of $X$ is
called \textit{ $\varepsilon$-absorbing} if $U$ contains the
$\varepsilon$-ball centered at $f(x)$ for each $x$ in $U$. $U$ is
\textit{absorbing} if it is $\varepsilon$-absorbing for some
$\varepsilon >0$. In an equivalent description, $f$ maps $U$ in a
uniform distance $\varepsilon$ into its interior.
\end{defn}

\smallskip

\begin{defn} Let $U$ be a nonempty absorbing subset of $X$. The closed set $$A=\cap_{n \geq 0}
\overline{f^{n}(U)}$$ is called the \textit{attractor-like set
determined by $U$}. In the compact case, $A$ is usally referred to
simply as the \textit{attractor} determined by $U$. The open set
$\cup_{n\geq 0}f^{-n}(U)$ is the set of all points whose
omega-limit sets are contained in $A$ and it is called the
\textit{basin of $A$ relative to $U$, denoted by B(A;U)}. It is
well known that, if $X$ is compact, $B(A;U)$ is independent of the
choice of absorbing open sets $U$.
\end{defn}

\smallskip

\smallskip

\begin{defn}A point $x$ in $X$ is  \textit{nonwandering} with
respect to $f$ if  any open neighborhood $U$ of $x$, there exists
a positive integer $N$ such that $f^{N}(U) \cap U \neq
\varnothing$. The set of all nonwandering points of $f$ is denoted
by  \textit{$NW(f)$}.
\end{defn}

\smallskip

The nonwandering set $NW(f)$ is closed and positively invariant.
If, in particular, $f$ is bijective, then $f(NW(f))=NW(f)$. In the
compact space, the nonwandering set is nonempty.

\smallskip

\begin{defn}
An \textit{$\varepsilon$-chain (= $\varepsilon$-pseudo-orbit)}
\textit{from} $x_0$ \textit{to} $x_n$ for $f$ is a sequence $x_0,
x_1,\cdots , x_n$ with the property that $d(f(x_i), x_{i+1}) <
\varepsilon $ for $0 \leq j < n.$ A point $x$ is a \textit{chain
recurrence point} if for every $\varepsilon >0$, there exists an
$\varepsilon$-pseudo-orbit with the starting and ending points at
$x$, that is, there exist $n\in\zz_+$ and $x_0,x_1,\cdots,x_n\in X$
such that $x=x_0,x_n$ and $d(f(x_i), x_{i+1}) < \varepsilon $ for
all $i=0,\cdots,n-1$. The \textit{chain recurrence set $CR(f)$} of
$f$ is the set of all chain recurrence points, i.e.,
$$
CR(f)= \{p \in X|\ \mbox{there exists an $\varepsilon $-chain from
$p$ to itself, for every $\varepsilon > 0$}\}.
$$ \end{defn}\smallskip
It is obvious that the chain recurrence set $CR(f)$ is closed. This
set was defined by C. C. Conley \cite{Co}. He showed that in case X
is a compact metric space, $CR(f)$ is globally determined by the
attractors and the basin of attractors. The following theorem is
Conley's characterization of the chain recurrence set for continuous
maps:

\smallskip

\begin{lem}\label{thm:co} ({\rm Conley's Theorem} \cite{Co})
Let $X$ be a compact metric space and $f:X\rightarrow X$ be
continuous. Then the chain recurrence set of $f$ is the complement
of the union of the set $B(A)-A$, as $A$ runs over the collection of
attractors of $f$; here $B(A)$ denotes the basin of attractor $A,$
that is,
\begin{eqnarray*}
X-CR(f)= \bigcup_{A: attractor} ( B(A)-A )
\end{eqnarray*}
\end{lem}

The volume-preserving condition obstructs the existence of
attractors, and hence the existence of non-chain recurrences, as
below.

\begin{prop}\label{thm:cptcase}
Let $M$ be a compact  manifold with a volume form $\omega$ and $f$
be a volume-preserving diffeomorphism on $M$. Then, $M$ is a chain
recurrence for $f$, i.e., every point of $M$ is a chain recurrence
for $f$.
\end{prop}

\begin{proof}
By the definition of the volume-preserving diffeomorphisms, there
does not exist an absorbing set, because any subset is of finite
measure. Thus, there does not exist an attractor. Therefore the
right side of the equation in Lemma \ref{thm:co} is empty. Hence
$M = CR(f)$, i.e., $M$ is chain recurrent for $f$.
\end{proof}

\smallskip

As one deals with non-compact manifolds, the definitions around
attractors should be corrected, which we do in the next section.
However, in the volume-preserving dynamics on non-compact manifolds
with a volume form $\omega$, in general, we do not guarantee the
parallel consequence of the theorem \ref{thm:cptcase} (Example
\ref{ex: noncpt}).

\smallskip


\begin{rk}
The following theorem exhibits that a compact manifold $M$ with a
volume form $\omega$ and a volume-preserving diffeomorphism $f$ is
nonwandering. It is easy to see that the nonwandering set is
contained in the chain recurrence set, that is, $NW(f) \subseteq
CR(f)$. Using these facts, we directly say that the above $M$ is
chain recurrent.
\end{rk}

\smallskip

\begin{thm}\label{th: nonwander}
Let $M$ be a compact manifold with a volume form $\omega$ and $f$
be a volume-preserving diffeomorphism on $M$. Then $M=NW(f)$,
i.e., every point of $M$ is nonwandering.
\end{thm}

\begin{proof}
Assume that the conclusion is contrary. Thus there exists a
wandering point $p$ of $M$. Then we have some neighborhood $U$ of
$p$ satisfying that $f^{n}(U) \cap U = \varnothing$ for all
positive integer $n$. In fact, the positive orbit of $U$ under $f$
is composed the disjoint union of each iteration of $U$ under $f$.
If the union is not disjoint, then there are positive integers
$n_{1}, n_{2} (n_{1} > n_{2})$ such that $f^{n_{1}}(U) \cap
f^{n_{2}}(U) \neq \varnothing$. Then we obtain that
$f^{n_{1}-n_{2}}(U) \cap U \neq \varnothing$, which is
contradiction. Therefore $M$ contains the disjoint union
$\bigcup_{n\geq 0} f^{n}(U)$ and so
$$
m(M)\geq \sum_{n\geq 0} m(f^{n}(U))= \sum_{n\geq 0} m(U)=\infty.
$$
Hence the compactness of $M$ implies a contrary. This completes
the proof.
\end{proof}

\smallskip

From Theorem \ref{th: nonwander}, in the volume-preserving
version, we can obtain a stronger concept, say the notion of nonwandering.
However, nevertheless, we mainly treat the notion of chain
recurrence in the rest of this section because we have useful
tools to develop arguments.

Moriyasu \cite{Mo} proved that the interior of the set of elements $f$
in $Diff^{1}(M)$ with the property that the restriction of $f$ to the
nonwandering set has the shadowing property, is included in
$\mathcal{F}(M)$. Here, $\mathcal{F}(M)$ is the set of elements $f
\in Diff^{1}(M)$ such that there exists a
$\mathcal{C}^{1}$-neighborhood $\mathcal{U}(f)$ of $f$ satisfying
that for every element $g$ of $\mathcal{U}(f)$, all periodic
points of $g$ are hyperbolic. It is inferred from \cite[pp.3]{WGW}
and \cite{We} that $CR(f)$ is $\mathcal{C}^1$-stably shadowable if
and only if $CR(f)$ is hyperbolic. Therefore, we obtain the
following theorem.

\smallskip

\begin{thm}\label{th: sh hyp}
Let $M$ be a compact manifold with a volume form $\omega$ and $f$ be
a volume-preserving diffeomorphism on $M$. Then, $M$ is
$\mathcal{C}^{1}$-stably shadowable if and only if $M$ is
hyperbolic.
\end{thm}

\begin{proof}
This follows from Theorem \ref{thm:cptcase} and the result in
\cite{Mo} (see \cite[pp.2--3]{WGW}).
\end{proof}

\smallskip

We use (and have used implicitly) the uniform metric on the space
$Diff^{1}_{\omega}(M)$ of volume-preserving diffeomorphisms, given
by $d_U(f,g):=\sup_{x \in M}d(f(x),g(x))$. A subset $R$ is called
\textit{residual} if it contains a countable intersection of open
and dense subsets. If there exists a residual set $R$ in
$Diff^{1}_{\omega}(M)$ such that any $g$ in $R$ possesses the same
property $P$, then we call the property $P$ \textit{generic}. Such a
residual set is called a \textit{generic set}.

Bonatti and Crovisier\cite{BC} proved that, generically, a
volume-preserving diffeomorphism is transitive in a compact,
connected manifold with a volume-preserving diffeomorphism.

\smallskip

\begin{thm} \label{vol-transitive}
\textrm{(Theorem 1.3. in \cite{BC})} Suppose that M is a compact
connected manifold with a volume form $\omega$. Then, there exists a
residual set $\mathcal{G}_{\omega}$ in $Diff^{1}_{\omega}(M)$
consisting of transitive diffeomorphisms.
\end{thm}

\smallskip

Let $f$ be a diffeomorphism on $M$. An invariant subset $A$ of $M$
is \textit{chain transitive} if for every two points $p, q$ in $A$
and for every $\varepsilon
>0$, there exists an $\varepsilon$-pseudo-orbit with the starting
point $p$ and the end point $q$, that is, a finite sequence $p=x_0,
x_1,\cdots , x_n=q$ such that $d(f(x_i), x_{i+1}) < \varepsilon$ for
all $i$. The sequence $\{x_0, x_1,\cdots , x_n \}$ is called an
\textit{$\varepsilon$-chain in $A$ connecting $p$ and $q$}.

\smallskip

\begin{thm}\label{thm:chain-transitive}
Let $M$ be a connected and compact manifold with a volume form
$\omega$, and $f$ be a volume-preserving diffeomorphism on $M$. Then
$M$ is chain transitive for $f$.
\end{thm}

\begin{proof}
Firstly we denote the transitive residual set in Theorem
\ref{vol-transitive} by $R$. For every $\varepsilon >0$, there
exists a positive number $\delta$ with $\delta < \varepsilon$  such
that
$$
\mid x-y  \mid < \delta \ \mbox{implies}  \ \mid f(x)-f(y)\mid <
\varepsilon.
$$
From the Baire category theorem, we can pick a volume-preserving
mapping $g$ in $R$ such that $d_U(f,g)=\sup_{x\in M}d(f(x),g(x))<
\frac{\delta}{2}$. By the definition of $R$, $g$ is transitive,
i.e., there is a point $p_{0}$ in $M$ whose orbit closure is the
whole manifold $M$. For every point $p, q$, since the orbit of
$p_{0}$ is dense, there are positive integers $m, n$ such that
$g^{m}(p_{0}) \in B(p, \frac{\delta}{2})$ and $g^{n}(p_{0}) \in B(q,
\frac{\delta}{2})$. We may assume $n > m$. Now we consider the
sequence $\{ p, f^{m+1}(p_{0}), f^{m+2}(p_{0}), \cdots,
f^{n-1}(p_{0}), q \}$. Thus, by the uniform continuity of $f$, the
above sequence is an $\varepsilon$-chain in $M$ connecting $p$ and
$q$. This completes the proof.
\end{proof}

\smallskip

\begin{rk}
It is clear that every point of the chain transitive set is a chain
recurrence. Thus, in the connected case, not using the theorem
\ref{thm:cptcase} we are also able to say that $M$ is chain
recurrent for $f$.
\end{rk}

\smallskip
Next we state the case of symplectic manifolds.
\smallskip

\begin{rk}
Let $M$ be a connected and compact manifold and $\omega$ be a
symplectic form defined on $M$ and $Symp^{1}_{\omega}(M)$ the set of
symplectic $\mathcal{C}^{1}$-diffeomorphisms of $M$. \cite{ABC}
shows that $\mathcal{C}^1$-generic symplectic diffeomorphisms of a
compact manifold are transitive, that is, the transitive symplectic
diffeomorphisms form a residual subset $\mathcal{G}$ of
$Symp^{1}_{\omega}(M)$. \end{rk}

The following theorem also generalizes a previous result in
\cite{BC} assuring the density of transitive diffeomorphisms in
$Diff^1_\omega(M)$.

\begin{thm} (\cite[Theorem 1]{ABC}) If
$(M,\omega)$ is a compact connected symplectic manifold. The set $G$
of transitive symplectic diffeomorphisms forms a dense
$G_{\delta}$-set in $Symp^{1}_{\omega}(M)$. Moreover, the subset $G'
\subseteq G$ of diffeomorphisms with the unique homoclinic classes
$M$, forms a dense $G_{\delta}$-set in
$Symp^{1}_{\omega}(M)$.\end{thm}

\smallskip

\begin{rk}
By the similar method of the proof of Theorem
\ref{thm:chain-transitive}, we can say that the compact connected
symplectic manifold is also chain transitive.
\end{rk}

\smallskip

Let $f$ be a diffeomorphism on a compact manifold $M$. Let $\Lambda
\subseteq M$ be a closed $f$-invariant set. We say the subset
$\Lambda$ of $M$ is \textit{locally maximal} in $U$ if there exists
a compact neighborhood $U$ of $\Lambda$ such that $\cap_{n \in
\mathbb{Z}}f^{n}(U) = \Lambda$.

\smallskip

\begin{defn} A subset $\Lambda$ of $M$ is \textit{$\mathcal{C}^{1}$-stably shadowable}
if $\Lambda$ is locally maximal in some compact neighborhood $U$
and  there exists a $\mathcal{C}^1$-neighborhood
$\mathcal{U}(f)$ of $f$ such that for any $g \in \mathcal{U}(f)$, $g|\Lambda_{g}$ has the shadowing property
where $\Lambda_{g} = \cap_{n \in \mathbb{Z}}g^{n}(U)$.
\end{defn}

\smallskip

The above set $\Lambda_{g}$ is called the \textit{continuation} of
$\Lambda$. Specially, if $\Lambda=M$, we call simply $f$ is
\textit{$\mathcal{C}^{1}$-stably shadowable}.

\smallskip

Let $p, q$ be points in $M$. Then,
for every $\varepsilon>0$,
we can consider an equivalence
relation $R$ on $CR(f)$ as follows. $pRq$ means there exist both an
$\varepsilon$-chain from $p$ to $q$ and an $\varepsilon$-chain from
$q$ to $p$. Then we call the equivalence classes the \textit{chain
recurrence classes} or the \textit{chain transitive components} of
$f$, simply call the \textit{chain components}. The components are
compact invariant sets and cannot be decomposed into two disjoint
compact invariant sets. Denote $C_{f}(p)$ the chain component of $f$
that contains $p$. It is natural that $C_{f}(p)$ is
$\mathcal{C}^{1}$-stably shadowable if there is a $\mathcal{C}^{1}$-
neighborhood $\mathcal{U}$ of $f$ such that for every $g \in
\mathcal{U}$, $C_{g}(p_g)$ has the shadowing property, where $p_g$
is the continuation of $p$.

\smallskip
The following theorem by Wen et.al. is useful.
\smallskip

\begin{thm}\label{thm:hyp} \cite{WGW}
Let $p$ be a hyperbolic periodic point of $f$.
If $C_{f}(p)$ is $\mathcal{C}^{1}$-stably shadowable, then
$C_{f}(p)$ is hyperbolic.
\end{thm}

\smallskip

\begin{thm}\label{thm:global hyp}
Let $M$ be a compact manifold with a volume form $\omega$ and $f$ be a
volume-preserving diffeomorphism on $M$. Assume that every
chain component has a hyperbolic periodic point. If every chain component
is $\mathcal{C}^{1}$-stably shadowable, then $M$ is hyperbolic.
Furthermore, the converse is also true.
\end{thm}

\begin{proof}
Since the number of connected components of $M$ is finite,  there
exists $n\in\zz_{>0}$ such that $f^n$ maps each component to itself.
I.e., the restriction of $f^n$ to any connected component is also a
volume-preserving diffeomorphism on the component. From Theorem
\ref{thm:chain-transitive}, the connected component of compact
manifold $M$ is chain transitive for $f^n$. By the invariance, the
connected component becomes a chain component for $f^{n}$. Note that
every chain component of $f^n$  is included in some chain component
of $f$. Thus we can also say that every chain component of $f$ is
some union of connected components of $M$. Then every chain
component of $f$ is clopen (=closed and open) in $M$. By the
assumption of this theorem and Theorem \ref{thm:hyp}, the chain
component is also hyperbolic with respect to $f$. Since $M$ is chain
recurrent, it consists of finite union of some clopen chain
components which are hyperbolic. Then we conclude the first
implication. The converse follows directly from the definitions.
\end{proof}

\smallskip

\begin{rk}
In general, the number of the chain components of a compact manifold
is infinite. However, from the proof of the above theorem, in
the volume-preserving case, we can say that the number of the
components of $M$ is just finite. More precisely, $M$ is composed of
finitely many chain transitive components and every component is
also composed of finitely many connected components.
\end{rk}

\smallskip


\section{Chain recurrences on Non-compact manifolds}\label{sec: chain2}

We start with some necessary notions which are different from those
for the compact manifolds. We assume that $(X, d)$ is a metric space
and $f:X \rightarrow X$ is a homeomorphism. We define
$$\mathcal{P}=\mbox{the set of $\rr^+$-valued continuous functions on $X$}.$$

\smallskip
\begin{defn}
A nonempty open subset $U$ of $X$ is \textit{weakly absorbing} for
$f$ if there exists a mapping $\varepsilon \in \mathcal{P}$ such
that the closed ball $B_{\varepsilon(f(x))}(f(x)) \subseteq U$ for
each $x\in U.$ When $U$ is weakly absorbing, the set $$A=\bigcap_{n
\geq 0} \overline{f^{n}(U)}$$ is the \textit{weak attractor
determined by $U$}. If $\varepsilon \in \mathcal{P}$, then $x_0.
x_1, \cdots, x_n$ is an \textit{$\varepsilon(x)$-chain} if
$d(f(x_j), x_{j+1})<\varepsilon(f(x_j))$ for $0 \leq j <n-1$. The
number $n$ is called the \textit{length} of the
$\varepsilon(x)$-chain. A point $p$ is \textit{strongly chain
recurrent} for $f$ if for every $\varepsilon \in \mathcal{P}$, there
exist an $\varepsilon(x)$-chain of length at least $1$ that begins
and ends at $p$. We denote by
$$CR^{+}(f)= \mbox{the set of all strong chain recurrence points of $f$ }. $$
\end{defn}
\smallskip

In this section, for brevity, an attractor always means a weak attractor and also an attractor block or an absorbing set means a weakly absorbing set.
It is easily proved that $U$ is weakly absorbing if and only
if $\overline{f(U)} \subseteq U.$

\begin{defn}
In the case of non-compact spaces,  we define the \textit{basin of
an attractor $A$ relative to $U$, B(A;U)} as the open set
$\cup_{n\geq 0}f^{-n}(U)$. Every point of $B(A;U)$ has the
omega-limit sets contained in $A$. When $X$ is a compact space,
$B(A;U)$ is independent of $U$ while it is not true for non-compact
manifolds. Therefore, we define the \textit{extended basin} $B(A)$
of $A$ by the union of the set $B(A;U)$ as $U$ runs over all
the absorbing sets that determine $A$.  \end{defn}

\smallskip

Let us get back to Proposition \ref{thm:cptcase} for a while. The
example below exhibits the failure of the theorem in the
volume-preserving dynamics over non-compact manifolds.

\begin{exam}\label{ex: noncpt}
Let $M=\mathbb{R}^2$ and $f : M \rightarrow M$ given by
$f(x,y)=(x+1,y)$. Let $\omega$ be a volume form (equivalently, a
symplectic form) by $\omega=dx\wedge dy$. Then, it is clear that $f$
preserves $\omega$. Let
$$ U_n=\{(x,y)\in M \ | \ y < \frac{-1}{x-n}, \ x <
n \} \cup \{(x,y)\in M \ |\ x \geq n \}.$$ Since $f(U_n)=U_{n+1}$,
we can easily check that $U_0$ is an attractor block for the
translation $f$ and $$A=\{(x,y) \ | \ y\leq 0 \}$$ is the attractor.
Whilst, no point of $M$ is a (strong) chain recurrence for $f$.
\end{exam}

The following theorem by Hurley is a generalized version of Conley's
theorem. In \cite{Hu1,Hu2}, Hurley investigates the extent to which
Conley's Theorem can be salvaged, and generalized Conley's theorem
from dynamical systems on compact metric spaces to those on locally
compact metric spaces. Thus, his theorem is applicable to an
arbitrary manifold.

\smallskip

\begin{thm}\label{thm:noncpt Hurley}\cite{Hu1,Hu2}
If $X$ is a locally compact metric space and $f:X\rightarrow X$ is
continuous, then the strong chain recurrence set $CR^{+}(f)$ of $f$
is the complement of the union of the set $B(A)-A$, as $A$ runs over
the collection of weak attractors of $f$. I.e.,
\begin{equation}\label{eq: Hurley}
X-CR^{+}(f)= \bigcup_{A: weak\ attractor} ( B(A)-A ).
\end{equation}
\end{thm}

\smallskip

\begin{rk}
A topological dynamical system $(X, f)$ is called $\textit{minimal}$ if
the orbit of every point $x \in X$ is dense in $X$. Let
$f:X\rightarrow X$ be a continuous map on a metric space $X$, and
$(X, f)$ be a minimal dynamical system. From the definition of the
minimality, $M$ is strongly chain recurrent for $f$, i.e., every
point of $M$ is strongly chain recurrent for $f$.
\end{rk}

\smallskip

Now we have a generalization of Propositon \ref{thm:cptcase} to
the case of compact attractors.
\begin{thm}
 Let $M$ be a manifold (not necessarily compact) with a volume form $\omega$,
and $f$ be a volume-preservig diffeomorphism on $M$. If every
attractor of $M$ is compact, then $M$ is strongly chain recurrent
for $f$.
\end{thm}

\begin{proof}
Assume that there exists an element $p$ of $M$ that is not
strongly chain recurrent with respect to $f$. By Hurley's theorem,
there is an attractor $A$ and a basin $B(A)$ of $A$ such that $p
\in B(A)-A$. Since $A$ is compact, we may assume that there is an
attractor block whose closure is compact. From now on, the
remaining part of the proof is parallel with the proof of
Proposition \ref{thm:cptcase} for the compact manifolds.
\end{proof}

\smallskip
The following proposition (and its corollary) shows the invariance of the attractors and the boundaries.
\begin{prop}\label{prop:orbit}
Let $f$ be a homeomorphism on a metric space $X$, $U$ be an
attractor block, and $A$ be a weak attractor depending on the
attractor block. If a point $x$ is in $U-A$, then the intersection
of the (positive) $f$-orbit of $x$ and the attractor is empty.
\end{prop}

\begin{proof}
Let $O^+_f(x)$ be the (positive) $f$-orbit of $x$. Suppose $O^+_{f}(x)
\cap A \neq \varnothing$, Then there exists a nonnegative integer
$k$ such that $f^{k}(x) \in A$, that is, $f^{k}(x) \in \cap_{n \geq
0} \overline{f^{n}(U)}$.
Note that, $f(\overline{U})=\overline{f(U)}$.
Thus $x \in f^{n-k}(U)$ for all $n \geq 1$
and so $x \in A$ by the shrinking property. This is a contradiction,
which completes the proof.
\end{proof}

\smallskip

When $f:X\to X$ is continuous, by the definition, it is easily shown
that an attractor is positively $f$-invariant. If $f$ is a
homeomorphism, an attractor $A$ is $f$-invariant, i.e., $f(A)=A$.
Indeed, if $f(A)\neq A$ then there is an element $x$ in $A - f(A)$.
From the definitions, $f^{-1}(x) \in U-A$, where $U$ is an
associated attractor block. Then, Proposition
\ref{prop:orbit}, we must meet a contradiction. Hence an attractor
is invariant.
\smallskip

\begin{cor}
Let $f$ be a homeomorphism on a locally compact metric space $M$.
Then the boundary of every weak attractor is positively $f$-invariant, that is,
$f(\partial A)\subseteq \partial A$ for every weak attractor $A$.
\end{cor}

\begin{proof}
Suppose the contrary of the conclusion. Then by the above statement,
we may assume that there exists a  boundary point $x$ satisfying $f(x)$ is in the interior of $A$.
From the local compactness, we can choose compact neighborhood $C$ of $f(x)$
such that $f^{-1}(C)$ is also a compact neighborhood of $x$.
Then, we are able to pick a point in $U-A$ where $U$ is an associated attractor block which corresponds to an interior point of the attractor.
By Proposition \ref{prop:orbit}, it is a contradiction.
\end{proof}

\smallskip


Now we embark on the main proposition and theorem for the chain
recurrences on the non-compact cases. The first proposition tells us
that the points near an attractor with bounded orbits form a measure
0 set, in the volume-preserving dynamics. Recall Example \ref{ex:
noncpt} which shows every orbit is unbounded.

\begin{prop}\label{prop: measure 0} Let $M$ be a manifold (not
necessarily compact) with a volume form $\omega$. Let $f$ be a
volume-preserving diffeomorphism on $M$. Let $A$ be any attractor
and $U$ be an associated attractor block.
Then, the complement in $U-A$ of the set of points $p\in U-A$ with
unbounded orbits is of measure 0.
\end{prop}

\proof Let $p\in U-A$ and $K\subset U-A$ be a compact neighborhood of $p$ with a finite measure $c>0$. Let us fix any point $x_0\in M$.
Let $B_l(x_0)$ be the closed ball of the radius $l$ centered at $x_0\in M$ (where $l\in \mathbb{Z}_+$).
Let us define
\begin{equation}
 K_r= \{q\in K| f^k(q) \notin B_r(x_0)\ \mbox{for some positive integer $k$} \}
\end{equation}

It suffices to show \begin{equation}\label{eq: L Kr} L=\bigcap_{r\in
\mathbb{Z}_+}K_r\end{equation} is of measure $c$, because $L$ is the
set of points of $K$ with unbounded orbits and $m(L)=m(K)$ implies
the statement of the proposition. Note that

\begin{eqnarray*}
K-K_r &=& \{q\in K| f^k(q) \in B_r(x_0)\ \mbox{for all $k\in \mathbb{Z}_+$} \}  \\
&=& \bigcap_{k\in \mathbb{Z}_+} \{q\in K| f^k(q) \in B_r(x_0)  \}
\end{eqnarray*}
and thus $K-K_r$ is measurable as $\{q\in K \ | \ f^k(q) \in B_r(x_0)  \}$
is measurable for each $k\in\zz_+$. Therefore, $L$ is measurable as well.

Let us observe that \begin{equation}\label{eq: 0}
m((f^k(U)-f^k(A))\cap B_r(x_0)) \to 0 \end{equation} as $k\to
\infty$. Indeed, Lebesgue's dominated convergence theorem assures it
from the following: \begin{itemize}\item[(a)] definition of
attractors (i.e., $\cap_k\overline{f^k(U)}=A$), \item[(b)] the $f$-invariance
of $A$, \item[(c)] $B_r(x_0)$ is of finite measure.
\end{itemize}  Note that $\{q\in K| f^k(q) \in B_r(x_0)  \}=f^{-k}(B_r(x_0))\cap
K$. Thus, we have \begin{eqnarray*} m(\{q\in K| f^k(q) \in B_r(x_0)
\}) &=& m(f^{-k}(B_r(x_0))\cap K) \\ &=& m( (B_r(x_0))\cap
f^k(K))\end{eqnarray*} where the latter equality is due to the
measure-preservation of $f$. Because of the inclusion $f^k(K)
\subset f^k(U)-f^k(A)$ and \eqref{eq: 0}, we obtain
$$m(\{q\in K| f^k(q) \in B_r(x_0)  \}) \to 0$$ as $k\to \infty$ by
Lebesgue's dominated convergence theorem. Therefore, for each $r$,
we have $m(K-K_r)=0$, equivalently, $m(K_r)=m(K)-m(K-K_r)=c$. By
applying Lebesgue's dominated convergence theorem to \eqref{eq: L
Kr}, we obtain $m(L)=c$, as
desired. \qed\\

For $p \in M$, we denote by $O^+_{f}(p)= \{f^{n}(p) \ | \ n \geq 0
\}$ and we define that $K^+(p)= \overline{O^+_{f}(p)}$. $M$ is said
to be \textit{Lagrange-stable} for $f$ if for $p \in M$, $K^+(p)$ is
compact. Since we are working on a metric space, the Lagrange
stability amounts to $O^+_f(p)$ is bounded.

\smallskip

\begin{thm}\label{th: LS CR} Let $M$ be a
manifold with a volume form $\omega$, and $f$ be a Lagrange-stable
volume-preserving diffeomorphism on $M$. Then, $M$ is strongly chain
recurrent for $f$, i.e., every point of $M$ is strongly chain
recurrent with respect to $f$.
\end{thm}

\begin{proof} This follows from Hurley's theorem (Theorem
\ref{thm:noncpt Hurley}) for locally compact spaces. To prove the
theorem, as was shown in Theorem 4.2, the nonexistence of attractors
should be guaranteed. On the contrary, suppose that a nonempty
attractor $A$ exists. By the above proposition, almost every point
of $U-A$ has an unbounded orbit, where $U$ is an attractor block of
$A$. This contradicts to our assumption of the Lagrange stability.
\end{proof}

\smallskip


\section{Nonshadowability and Recurrence}\label{sec: rec}

In this section, we show the theorems about the nonshadowablity and
recurrence on compact manifolds.

\begin{lem}\label{lemma:stable}
Let $p\in M$ be a hyperbolic periodic point. If for $x\in W^s
(p,f),$ the positive orbit starting from $x$ is $\epsilon$-shadowed
by an actual orbit for all sufficiently small $\epsilon>0,$ then the
actual orbit is in the stable manifold.
\end{lem}

\begin{proof} The lemma is proven in a small neighborhood near $p$;
the problem is linearized as follows. Let $f:\rr^n\to\rr^n$ be a
hyperbolic linear map. Let $x$ lie in the direct sum $\Sigma$ of
eigenspaces with eigenvalues of absolute values $<1$. Because the
positive orbit of any point $y$ outside $\Sigma$ diverges, the
positive orbit of $x$ cannot be $\epsilon$-shadowed by the positive
orbit of $y$.
\end{proof}

\smallskip

With the above lemma, we can prove the following :

\begin{thm}\label{th: rec nonsh}
Let $M$ be a compact manifold with a volume form $\omega$ and $f : M
\rightarrow M$ be a volume-preserving diffeomorphism. If $W^u (p, f)$ has a recurrence point
and there is no homoclinic point of $p,$ then $f$ is nonshadowable.
\end{thm}

\begin{proof}
Let $q$ be a recurrence point in $W^u (p, f).$ Then, there exists a
point $x_0$ in $W^s (p, f)$ such that $x_0$ is in the omega-limit
set of $q$ (see \cite[Proof of Theorem 1]{Xi1}.) To show
nonshadowability, it suffices to find a $\delta$-pseudo-orbit for
each $\delta > 0$ which is not $\epsilon$-shadowed for some
$\epsilon
>0.$ First, select a point $f^{n_0}(q)$ in the positive orbit of $q$
such that the distance between $f^{n_0}(q)$ and $x_0$ is less than
$\delta.$ Then, we consider the following $\delta$-pseudo-orbit
$(y_i)$ defined as
\begin{enumerate}
  \item for $i < n_0 -1, y_i = f^i (q),$
  \item for $i \geq n_0, y_i = f^{i-n_0} (x_0).$
\end{enumerate}
By applying Lemma \ref{lemma:stable} respectively to the positive
orbit and the negative orbit of $(y_i),$ it is easily obtained that
there exists a positive number $\epsilon$ such that the
$\delta$-pseudo-orbit $(y_i)$ cannot be $\epsilon$-shadowed by any
actual orbit because $p$ has no homoclinic point.
\end{proof}

\smallskip

Because by \cite{Xi1},  a hyperbolic periodic point has a homoclinic
point for a generic volume-preserving diffeomorphism, the above
theorem seems uninteresting. But, there is an interesting example in
\cite[p. 655]{Ol1} which admits a recurrence point but has no
homoclinic point. Hence, our theorem is not vacuous and guarantees
that the example is nonshadowable\cite{Ol3}.

With the similar arguments, we can obtain a result for the case when
there is a homoclinic point under the shadowing property.

\smallskip

\begin{thm}
Let $M$ be a compact manifold with a  volume form $\omega$ and $f :
M \rightarrow M$ be a volume-preserving diffeomorphism. Let $p$ be a hyperbolic periodic point. Assume that
$f$ has the shadowing property. If there is a recurrence point in
$W^u (p, f),$ then the recurrence point lies in the limit set
of homoclinic points of $p.$
\end{thm}

\begin{proof}
Let $q$ be the recurrence point. We need to find homoclinic points
converging to $q.$ For each $n\in\zz_+$, we can construct a
$1/n$-pseudo-orbit by the exactly same method in the proof of the
above theorem. Then, by the shadowability of $f$, the pseudo-orbit
is $\epsilon$-shadowed by an actual orbit for every sufficiently
small $\epsilon>0.$ Especially, we may assume that $\epsilon$ is
less than $1/n.$ By applying Lemma \ref{lemma:stable} twice as in
the proof of Theorem \ref{th: rec nonsh}, the actual orbits are all
homoclinic points. Hence, we can find a homoclinic point in the ball
centered at $q$ of radius $1/n$.
\end{proof}

By this result, we can obtain the following directly.

\begin{cor}
Let $M$ be a compact manifold with a  volume form $\omega$ and $f :
M \rightarrow M$ be a volume-preserving diffeomorphism. Let $p$ be a hyperbolic periodic point. If $f$ has
the shadowing property and the recurrence points are
dense in $W^u (p, f),$ then the homoclinic points are also dense in
$W^u (p, f).$ \qed
\end{cor}

\smallskip

In \cite{Ol2}, it is proved that the recurrence points of a generic
volume-preserving diffeomorphism on a compact manifold $M$ with
volume form $\omega$, are dense in an unstable manifold. Xia
\cite{Xi1} also showed generically hyperbolic periodic point $p$ of
$f$ has a homoclinic point, and moreover, the homoclinic points of
$p$ are dense in both the stable manifold and the unstable manifold
of $p$. By the corollary above, without genericity, we obtain a connection between the
recurrence points and homoclinic points under  the shadowability
condition.

\smallskip



\end{document}